\documentclass{elsart}

\usepackage{amssymb,amstext}
\usepackage[latin1]{inputenc}
\journal{\small IJNT}

\font\mini=cmr10 at 7pt

\def\q{\mathfrak q}

\def\dsp{\displaystyle}
\def\gal{\hbox{\rm Gal}}
\def\cl{\hbox{\rm Cl}}
\def\pcl{\hbox{\mini Cl}}

\def\hom{\hbox{\rm Hom}}
\def\Res{\hbox{\rm Res}\,}

\def\ind{\hbox{\rm Ind}}
\def\hom{\hbox{\rm Hom}}
\def\ker{\hbox{\rm Ker}}
\def\frob{\hbox{\rm Frob}}

\def\id{\hbox{\rm id}}
\def\ver{\hbox{\rm Ver}}
\def\Im{\hbox{\rm Im}}
\def\ov{\overline}
\def\No{{\rm N}}

\def\plus{\ds\mathop{\raise 2.0pt \hbox{$\bigoplus $}}\limits}
\def\mult{\ds\mathop{\raise 2.0pt \hbox{$\bigotimes$}}\limits}
\def\prd{ \ds\mathop{\raise 2.0pt \hbox{$  \prod   $}}\limits}
\def\Cap{ \ds\mathop{\raise 2.0pt \hbox{$\bigcap   $}}\limits}
\def\Cup{ \ds\mathop{\raise 2.0pt \hbox{$\bigcup   $}}\limits}
\def\sm{  \ds\mathop{\raise 2.0pt \hbox{$  \sum    $}}\limits}

\def \tensorZ {\otimes{\raise -0.8pt \hbox{\!\!$_{_{\zz}}$}}}
\def \tensorZp {\otimes{\raise -0.8pt \hbox{\!\!$_{_{\zz_{\!p}}}$}}}
\def \tensorQ {\otimes{\raise -0.8pt \hbox{\!\!$_{_{\qq}}$}}}
\def \tensorQp{\otimes{\raise -0.8pt \hbox{\!\!$_{_{\qq_{\!p}}}$}}}
\def \tensorCp {\otimes{\raise -0.8pt \hbox{\!\!$_{_{\cc_{\!p}}}$}}}
\def \tensorK {\otimes{\raise -0.8pt \hbox{\!\!$_{_{K}}$}}}

\def\fin{\vbox{\hrule\hbox to 7.2pt{\vrule height 7pt\hfil\vrule}\hrule}}

\def\semi{\mbox{$\times$}\mkern -4.9mu {\raise 1.2pt \hbox{\tiny $\vert$}}}
\def\lien{\mathrel{\mkern-4mu}}
\def\too{\relbar\lien\rightarrow}
\def\tooo{\relbar\lien\relbar\lien\too}

\let\ov=\overline

\begin{document}
\begin{frontmatter}

\title{Principalization of ideals in abelian extensions of number
fields}

\author{S\'ebastien {\sc Bosca}}

\address{Universit\'e de Bordeaux, Institut de Math\'ematiques de
Bordeaux, 351 cours de la
Lib\'eration, 33405 TALENCE Cedex, France}

\medskip
{\footnotesize With an Appendix by} Georges {\sc Gras}\footnote{\em Villa la Gardette, Chemin Ch\^ateau Gagni\`ere,
38520 Le Bourg d'Oisans, France} {\footnotesize and} Jean-Fran\c cois {\sc Jaulent}
\footnote{\em Universit\'e de Bordeaux, Institut de Math\'ematiques de
Bordeaux, 351 cours de la Lib\'eration, 33405 TALENCE Cedex, France.}

\begin{abstract}
We give a self-contained proof of a general conjecture of  G. Gras on
principalization of ideals in abelian extensions of a given field $L$, yet
solved by M.
Kurihara in the case of totally real extensions $L$ of the rational field
$\mathbb
Q$.\par
More precisely, for any given extension $L/K$ of number fields, in
which at least one infinite place of $K$ is totally split, and for
any ideal class $c_L^{}$ of $L$, we build a finite abelian extension
$F/K$, in which all infinite places  are totally split, such that
$c_L^{}$ principalizes in the compositum $M=LF$.
\end{abstract}
\end{frontmatter}

\medskip
\centerline{A -- INTRODUCTION}

\bigskip
When $M/L$ is an abelian extension of number fields, the problem of
knowing which ideals of $L$ are principal in $M$ is difficult, even if
$M/L$ is cyclic. Class field theory gives partial answers. For
instance, the Artin-Furtw\"{a}ngler theorem states that when $M=H_L$ is
the Hilbert class field of $L$, all ideals of $L$ are principal in $M$; but
the other cases are more mysterious.

When $M/L$ is cyclic, we still have no general answer but the problem
is easier. The kernel of the natural map $j:\cl_L^{} \rightarrow \cl_M^{}$
is
partly known by this way,
as explained below: at first, cohomology of
cyclic groups says that
this kernel is a part of $\hat H^1(\mathcal G,E_M^{})$, where $\mathcal
G=\gal(M/L)$ and
$E_M^{}$ is the group of units of $M$; $\hat H^1(\mathcal G,E_M^{})$ itself
is  not well
known but its order can be deduced from the order of $\hat
H^0(\mathcal  G,E_M^{})$
using the Herbrand quotient; after that, the order of $\hat
H^0(\mathcal G,E_M^{})$
depends partly from the natural map $E_L^{}\rightarrow U_S^{}$, where
$U_S^{}$ is the subgroup of the unit id\`eles of $L$ which is the product of
the groups of
local units at the places ramified in $M/L$ ($S$ is the set of
such ramified places); finally, the map $E_L^{}\rightarrow U_S^{}$ depends
on the Frobenius' of the primes of $S$ in the extension
$L[\mu_h^{}, \, E_L^{1/h}]/L$
(where $h=[M:L]$ and where $\mu_h^{}$ is the
group of $h$th roots of
unity). Obviously, this makes sense only if the primes of $S$
 are prime to the degree of the extension.

However, even in this cyclic case, $\ker(j)$ is not completely
given by these Frobenius', so that knowing such Frobenius', we cannot
get  really more than a
minoration of the order of $\ker(j)$.

This article deals with such minoration techniques, which allow to
prove, for instance, this easy fact: if a cyclic extension $M/L$ is
ramified at only one finite place $\q$ prime to $[M:L]$, then,
as soon as the ramification index $e_\q$ is large enough (precisely,
when $|\cl_L^{}|$ divides $e_\q$), the capitulation kernel $\ker j$
contains at least the class of $\q$ in $\cl_L^{}$. This can be easily
established by  studying $\hat H^1(\mathcal G,E_M^{})$ and may be seen
as a particular  case of our main theorem.

On the other hand, this article states a result which proves a
conjecture of Georges Gras and generalizes a result of Masato Kurihara (see
[1] and [2]), so that the theorem exposed here is not only an abstract
minoration of $\ker j$, using cohomology of cyclic groups,
Minkowski--Herbrand theorem on
units of number fields, class field theory, and Kummer
duality. Note also that here, we use new asymptotic methods (``take
$n$ large enough'', where $n$ is related to the degree $[F:K]$ in a suitable
manner).

Note finally that in the theorem below the hypothesis ``at least one
infinite place of $K$ is
totally split in $L/K$'' is necessary: in [1],
Georges Gras gives examples of extensions $L/K$ with ideals which do not
principalize in the compositum $LK^{ab}$,
where $K^{ab}$ is the maximal abelian extension of $K$.

\bigskip

\centerline{B -- MAIN THEOREM AND COROLLARIES}\medskip

\noindent {\bf Main Theorem.} {\em Let $L/K$ be a finite extension of
number fields in which at least one infinite place of $K$ totally splits.
There exists a
finite abelian extension $F/K$, which is totally split at all
infinite places, such that every ideal  of $K$ principalizes in the
compositum $M=LF$.}

\noindent{\bf Corollary 1} ($K=\mathbb Q$). {\em If $L$ is a totally
real  number field, any ideal
of $L$  principalizes in a real cyclotomic extension of $L$ (i.e.,
in the compositum of $L$ with a real subfield of a suitable
cyclotomic extension of $\mathbb Q$).}

Corollary 1 was proved by Masato Kurihara in [2].

\noindent {\em Notation:} In the following, $K^{ab}$ denotes the
maximal abelian extension of $K$ and  $K^{ab}_+$
its maximal totally real subextension  (i.e., the subextension of
$K^{ab}/K$ fixed under the decomposition groups of the infinite
places of $K$).

\noindent {\em Definition:} Let us  say that a number field $N$  is
principal when $\cl_N^{}$ is
trivial. Then:

\noindent{\bf Corollary 2.} {\em (i) Let $K$ be a number field with at least
one
complex place. Then any field containing $K^{ab}$ is
principal.\vspace{-0.35cm}

 (ii) Let $K$ be  a totally real number field. Any field
containing $K^{ab}_+$ and
whose Galois closure has at least one real place is principal.}

\noindent{\bf Corollary 3.} {\em (i) Any totally real field containing
$\mathbb Q^{ab}_+$ is
principal. \vspace{-0.35cm}

(ii)  Any field containing $\mathbb Q(i)^{ab}$ is principal.}

Corollary 3, for $\mathbb Q(i)$
or any imaginary quadratic field in (ii), was proved
by Masato Kurihara (see [2], theorem 1.1 p. 35 and theorem A.1 p. 46).

\noindent {\bf Corollary 4.} {\em
(i) Let $K$ be a number field with at least one
complex place. Then $K^{ab }$ is principal.\vspace{-0.35cm}

(ii) Let $K$ be a totally real number field. Any field
containing $K^{ab}_+$ which is contained in one
of the subfields of $K^{ab}$ fixed by a complex conjugation is
principal.}

Corollary 4 proves a conjecture of Georges Gras, Conjecture (0.5)
p. 405 of~[1].

\bigskip

\centerline{C -- PROOFS}

{\large \bf I. Proof of the main theorem : preliminaries}

To prove the theorem, we fixe an ideal $\mathfrak a_L$ of $L$, and we shall
build a finite abelian
extension $F$ of $K$, which is totally split at all infinite
places, cyclic in most cases, such that $\mathfrak a_L$ principalizes in the
compositum $LF$. Obviously, any ideal $\mathfrak b_L$ of $L$ with the same
class
in $\cl_L^{}$ will become principal in $LF$ as well, so that it is enough to
fix the class $c_L^{}$ of $\mathfrak a_L$ in $\cl_L^{}$.
Now, if $(c_i)_i$ is a finite generating system of $\cl_L^{}$, we will
obtain a corresponding set $(F_i)_i$ of
extensions, and every ideal of $L$ will principalize in $LF$, where
$F$ is the compositum of the $(F_i)_i$; this will prove the theorem.

So, in the following, we fixe $c_L^{}$ in $\cl_L^{}$, and must find $F$. As
the class group of $L$ is the direct sum of its $p$-parts, for all prime
numbers $p$, one can suppose that the order of $c_L^{}$ in $\cl_L^{}$ is a
power of
a prime $p$.
So, $c_L^{}$ and $p$ are fixed; $\cl_L^{}$ is now the $p$-part of
the class group of $K$, and $H_L$ the maximal $p$-extension contained
in the Hilbert class field of $L$.

{\bf(1) One can suppose $c_L^{}\in\cl_L^{p^a}$ for an arbitrary given
integer $a$:}

\noindent {\em Definition:} Let us call {\it abelian compositum} of the extension $L/K$
any extension $N=LF$, where $F$ is a finite abelian extension of $K$,
totally split at all infinite places of $K$.

One fixes an integer $a$; in case $c_L^{}\notin\cl_L^{p^a}$, one will
build an abelian compositum $L'$ of $L/K$, such that the extended class
$c_{L'}^{}=j(c_L^{})$ satisfies
$c_{L'}^{}\in\cl_{L'}^{p^a}$; so, if $N'$ is an abelian compositum of
$L'/K$ in which $c_{L'}^{}$ is principal,
it is as well an abelian compositum of $L/K$, so that one can
legitimately replace $L$ by $L'$, in which case one has
$c_{L'}^{}\in\cl_{L'}^{p^a}$.

Let's build such a $L'=LF_0$ as follows: let $\q$ be a prime
of $L$ satisfying the three following conditions:\vspace{-0.1cm}

(i) $\q$ totally splits in $L/\mathbb Q$;\vspace{-0.2cm}

(ii) $\frob(\q, L[\mu_{2p^a}]/L)=\id$;\vspace{-0.2cm}

(iii) $\frob(\q,  H_L/L)=c_L^{}$.\vspace{-0.1cm}

If such a $\q$ exists with say $\q | q$, the first two conditions
imply the existence of a (cyclic) subfield $F'_0$ of
$\mathbb Q(\mu_q)$, with degree $[F'_0:\mathbb Q]=p^a$, which is
totally ramified at the prime $q$;
the first condition implies that the compositums $F_0=KF'_0$ and
$L'=LF_0$
have again a degree $p^a$ over $K$ and $L$, respectively. Now $L'/L$ is
totally ramified at $\q$, say $\q = \q'\,{}^{p^a}$, for a prime $\q'$
in $L'$; so, according to the third condition, the extended class
$c_{L'}^{}=j(c_L^{})$ satisfies:
$$c_{L'}^{}=\overline\q=\overline{\q'}\,{}^{p^a}\in\cl_{L'}^{p^a}$$
\noindent as expected.

Now we only have to verify the existence of such a prime $\q$. The
three
conditions defining $\q$ all depend on $\frob(q,\tilde
H_L[\mu_{2p^a}]/\mathbb Q)$, where $\tilde H_L$ is the Galois closure
of $H_L$ over $\mathbb Q$; the first two conditions are equivalent to
the fact that
this Frobenius is in the subgroup $\gal(\tilde
H_L[\mu_{2p^a}]/L[\mu_{2p^a}])$; so they are compatible with the last
condition for any $c_L^{}$ in $\cl_L$, if and only if one has:
$L[\mu_{2p^a}]\cap  H_L=L$; if this
is right, the \v{C}ebotarev theorem states there are infinitely many
$\q$ satisfying the conditions. When this is wrong, we replace
$L$ by $L''=L[\mu_{2p^a}]\cap H_L$ which verifies
$L''[\mu_{2p^a}]\cap
H_{L''}=L''$.

As explained above, this last replacement is legitimate in case
$L''=L[\mu_{2p^a}]\cap H_L$ is an abelian compositum of $L/K$. In
fact, one has $L''=LF$ where $F$ is contained in the maximal
$p$-subfield of
$K[\mu_{2p^a}]$, which is clearly abelian and finite but in which
infinite
places are maybe not totally split for $p=2$.

So, for $p=2$, we shall complete the
proof, and we take in this particular case $L''=LK[\mu_{2^b}]_+$,
where the symbol $+$ denotes the maximal $\infty$-split subextension
over $K$, and where $b$ is an integer, which is choosen large enough
so that $L''$ contains $H_L\cap LK[\mu_{2^{a+1}}]_+$ and
$L''/L$ has degree at least 2.

Suppose we have found a prime $\q''$ of $L''$ satisfying the
following conditions:\vspace{-0.1cm}

(i) $\q''$ totally splits in $L''/\mathbb Q$;\vspace{-0.2cm}

(ii) $\frob(\q'', L''[\mu_{2^{a+1}}]/L'')=\id$;\vspace{-0.2cm}

(iii) $ \frob(\q'',  H_{L''}/L'')=c_{L''}^{}$,\vspace{-0.1cm}

where $c_{L''}^{}$ is extended from $c_L^{}$. So, let $F'_0$ be the totally
real subfield of $\mathbb Q[\mu_q]$ of
degree $2^a$ over $\mathbb Q$, thus $F_0=K[\mu_{2^b}]_+F'_0$ and
$L'''= LF_0= L''F'_0$  (which is an abelian compositum of $L/K$). If
$\q'''$ denotes the unique prime of $L'''$ above $\q''$, the extended
class $c'''_{L'''}$ in $\cl_{L'''}$ satisfies the expected condition:
$$
c_{L'''}^{}=\overline{\q''}=\overline{\q'''}\,{}^{2^a}\in\cl_{L'''}^{2^a}.
$$
So to conclude we only have to prove the existence of such a prime
$\q''$
verifying the threee conditions above. But this existence follows
from the \v{C}ebotarev theorem as soon as
the image of $c_L^{}\in\gal(H_{L''}/L'')$ in $\gal(H_{L''}\cap
L''[\mu_{2^{a+1}}]/L'')$ is trivial. To check this last point, let us
observe that in the class field description the extension of ideal
classes $j$ corresponds to the transfert map $\ver$.
Here $L''$ contains $H_L\cap LK[\mu_{2^{a+1}}]_+$, so $H_{L''}\cap
L''[\mu_{2^{a+1}}]=L'''$ is either $L''$ or $L''[i]$, and the image of
$j_{L''/L}(c_L^{})$ in $\gal(L'''/L'')$ is trivial, since one 
has:\vspace{+0.15cm}

\centerline{$\ver^{}_{B/A\rightarrow B'/A'}(\sigma)=\sigma^{[A':A]}$}\vspace{-0.15cm}

when $B'/A$ is abelian, so:\vspace{-0.15cm}
$$\Res_{L'''}(\ver^{}_{H_L/L\rightarrow
H_{L''}/L''}(c_L^{}))=\ver^{}_{H_L/L\rightarrow
L'''/L''}(c_L^{})=c_L^{[L'':L]}=\id.$$


{\bf(2) One can suppose that $ L/K$ is Galois:}

Let $\tilde L$ denotes the Galois closure of $L$ over $K$. Imagine the
theorem is proved for $\tilde L/K$ (in which at least
one infinite place totally splits as in $L/K$). Hence, there exists
an  abelian
compositum $\tilde LF$ of $\tilde L/K$ such that every ideal of
$\tilde L$ principalize in $\tilde LF$; so, $c_L^{}$ principalizes in
$\tilde LF$ but maybe does not in $LF$ and we have to study this case.

\smallskip
(a) When $c_L^{}$ is norm in $\tilde L/L$, say $c_L^{}=\No_{\tilde
L/L}(\tilde c_{\tilde L}^{})$,
the class  $\tilde c_{\tilde L}^{}\in \cl_{\tilde L}^{}$ principalizes in
$\tilde LF$,
say $j_{\tilde LF/\tilde L}(\tilde c_{\tilde L}^{})=1$ in $\cl_{\tilde
LF}^{}$; so we obtain:
$c_L^{}=\No_{\tilde L\cap LF/L}(c'_{\tilde L\cap LF})$ with $c'_{\tilde
L\cap LF}=(\No_{\tilde L/(\tilde L\cap LF)}(\tilde c_{\tilde L}^{}))$ and
the class $c'_{\tilde L\cap LF}$, which satisfies
$j_{LF/\tilde L \cap LF}\circ \No_{\tilde L/(\tilde L\cap LF)}(\tilde
c_{\tilde L}^{})=\No_{\tilde LF/LF}\circ j_{\tilde L F/\tilde L}(\tilde
c_{\tilde L}^{}) =1$, principalizes in $FL$;
 and so is $c_L^{}$.

 \smallskip
(b) When $c_L^{}$ is not a norm in $\tilde L/L$, maybe $c_L^{}$ is not
principal
in $LF$. But $\No_{\tilde L/L}(\cl_{\tilde L})$ contains
$\cl_L^{[\tilde L:L]}=\cl_L^{p^a}$, where $p^a$ is the largest power
of $p$ dividing $[\tilde L:L]$. According to Section {\bf (1)} we replace
$L$ by
$L'=LF_0$ such that $c_L^{}\in\cl_{L'}^{p^a}$. Since $[\tilde
L':L']=[\tilde LF_0:LF_0]$ divides
$[\tilde L:L]$, $c_L^{}$ is norm in $\tilde L'/L'$ and (a) applies.


\medskip
{\large \bf II. Proof of the theorem : building the extension $F$}

\noindent{\bf (3) The method and a first condition about the prime
$\q$: }

By now we suppose $L/K$ Galois. The prime $p$ and the class $c_L^{}$ are
fixed and we must
build an abelian compositum $LF$ of $L/K$ such that $c_L^{}$
principalizes
in $LF$. For convenience, we choose $F/K$ as a cyclic $p$-extension,
ramified at only
one finite place $\q$ of $K$, and whose  ramification index is
$e_\q(F/K)=p^n$ for a given integer $n$. We will see that with many
conditions about $\q$, when $n$ is large enough, $c_L^{}$ is principal in
$LF$, or in $L'F$, where $L'$ is a convenient abelian compositum of
$L/K$. At the end of the proof, in Section {\bf (6)}, we will study the
existence of such $\q$
verifying all conditions.

Now for a given integer $n$ and a given prime $\q$ of $K$, we wonder
if $c_L^{}$ is principal
in $M=LF$, where $F$ is a cyclic $p$-extension of $K$ with
$e_\q(F/K)=p^n$, unramified but at $\q$ and in which all
infinite places are totally split. The first question is the
existence  of such an
extension $F/K$ and class field theory gives the answer as follows.

Indeed, $N$ being the maximal abelian extension of $K$ unramified
but  at $\q$
and $\infty$-split, class field theory describes the Galois group
$\gal(N/K)$ from the
quotient
$$J_K\Big /K^\times.\prod_{v|\infty}K_v^\times.\prod_{\q'\neq
\q}U_{\q'}\ ,$$
where $J_K$ is the id\`ele group of $K$ and $U_{\q'}$ the
subgroup of local units of the completion $K_{\q'}$  of $K$ at the
place $\q'$. The inertia subgroup of
$\q$ in $N/K$ is, according to class field theory, isomorphic to the
quotient
$$U_\q \,\Big / \,
\overline{U_\q\cap(K^\times.\prod_{v \vert \infty}K_v^\times.\prod_{\q'\neq
\q}U_{\q'})}=U_\q\Big / \overline{E_K}\ ,$$
where $E_K$ is the group of global units in $K$ and the overlining
means closure in $U_\q$ of the diagonal embedding. Of course if $F$ exists,
it is contained in
the maximal
$p$-extension of $N$, whose ramification
subgroup is the $p$-part of $U_\q\Big / \overline{E_K}$, denoted
$(U_\q\Big / \overline{E_K})_p$.

We suppose now:\vspace{-0.2cm}

$\bullet \ \ \q\nmid p$.\vspace{-0.2cm}

So, if $F$ exists, $p^n$ divides $|(U_\q\Big / \overline{E_K})_p|$,
that is, under the assumption $\q\nmid p$:\vspace{-0.35cm}

$\bullet \ \mu_{p^n}\subset K_\q^\times\,$,\vspace{-0.35cm}

and\vspace{-0.35cm}

$\bullet \ E_K\subset U_\q^{p^n}\, .$\footnote{The canonical embedding of
$E_K$
in $K_\q^\times$
must be contained in $U_\q^{p^n}$ since  $(U_\q)_p$ is here a
cyclic group.}

These necessary conditions are enough to ensure the existence of $F$,
according to an obvious lemma, which states: if $A$ is an abelian
finite group and $C$ is a cyclic subgroup
of $A$ of which $p^n$ divides the order, then there exists a cyclic
quotient of $A$ in which the image of $C$ has order $p^n$.

Now we suppose that $\q \nmid p$ verifies the two above conditions
and the additional assumption:

\noindent $\bullet$ $\q$ is unramified in $L/K$.\footnote{In fact, in
the sequel Bosca will suppose that $\q$ is totally split in $L/K$.}

So $F$ exists;  all primes $\q_L^{} \mid\q$  are ramified in $M/L=LF/L$
with
the same index $e_{\q}=p^n$; and we have $[M:L]=p^{n+d}$ for some
positive integer $d$. \bigskip


{\bf (4) Obtaining a big cohomology group $ \hat H^0(\mathcal
G,E_M)$:}

Let $\mathcal G$ denotes the Galois group $\gal(M/L)$, which is
cyclic with  order
$p^{n+d}$, and $G=\gal(L/K)$.

According to the Minkowski--Herbrand theorem, the character of
 the representation
$\mathbb Q\otimes_\mathbb Z E_L$ of  $G=\gal(L/K)$,  given by the group of
global units 
$E_L$, is:
$$
\chi(E_L)=\sum_{v\mid \infty}\ind_{D_v}^G1_{D_v} -1_G\ ,
$$
where $1_G$ is the trivial character of $G$ and $1_{D_v}$ the trivial
character of the decomposition subgroup $D_v$.

Since at least one infinite place is totally split in the extension
$L/K$, the character of
$\mathbb Q\otimes_\mathbb Z (E_L/\mu_L E_K)$ satisfies
$$
\chi(E_L/\mu_L E_K)\ge\chi(\mathbb Z[G])-1\ ,
$$
and we can deduce from this the existence of a map:\footnote{See the
details in the Appendix.}
$$\varphi:E_L/\mu_L.E_K\longrightarrow\mathbb Z$$
such that, in $\hom(E_L/\mu_L.E_K,\mathbb Z)$, $\varphi$ generates a
$\mathbb Z[G]$-submodule whose character is $\chi(\mathbb
Z[G])-1$.

On the other hand, since $L[\mu_{p^n},E_L^{1/p^n}]/L[\mu_{p^n}]$ is a
Kummer extension, if $\q_L^{}$ is one of the primes of $L$
dividing $\q$, the Frobenius automorphism $$\sigma=\frob(\q_L^{},
L[\mu_{p^n},E_L^{1/p^n}]/L)$$ corresponds, in the Kummer duality,
to the map:
$$  u \mapsto ( u^{1/p^n})_{}^{(\sigma -1)}\in
\mu_{p^n}^{}, \ \ \hbox{for all $u\in E_L$}, $$
and we impose the new condition:

$\bullet\ \ \hbox{this map $\sigma$ coincides with
$\lambda_n:E_L\longrightarrow
E_L/\mu_L.E_K\stackrel{\varphi}{\longrightarrow}\mathbb
Z\longrightarrow\mu_{p^n}$}$,

\noindent
where the left map is the natural one and the right one is surjective
(we must choose a primitive $p^n$th root of unity for the right
map, but this choice does not change the Frobenius defined up to
conjugation: changing the choice of the root of unity is the same that
changing
the choice of a prime $\q' \mid \q$ in $L[\mu_{p^n}]$).

So, the property of $\varphi$ leads to the following facts:

Let $\phi$ denotes the map (see the Appendix):\vspace{-0.35cm}

\begin{eqnarray*}
\hspace{2cm} E_L &\  \longrightarrow & R_G:=
\{\dsp \sum_g\alpha_g \, g\in\mathbb Z[G]\  \mid \sum_g\alpha_g=0\} \\
u\  & \longmapsto & \ \ \dsp\sum_g\varphi(g^{-1}(u))\;g\  ;\end{eqnarray*}
$\Im(\phi)$ has finite
index, say $r$. So, $p^\delta$ being the maximal power of $p$ dividing
$r$, one has:
$$|\{E_L/ \{u\in E_L \  \vert\,  \forall g\in
G,\lambda_n(g(u))=1\}|\ge|R_G/R_G^{p^n}|/p^\delta=
p^{n(|G|-1)-\delta}\ .$$
On the arithmetical side, the ramification indices in $M/L$ of all
primes  $\q_L^{} \mid \q$ of
$L$ are all equal to $p^n$; so, with $\q_M^{} | \q_L^{} | \q$ in
$M/L/K$, one has\,{\footnote{Using the fact that the global norm is
the product of the corresponding local norms.}}:
$$\No_{M/L}(E_M)\subset
\No_{M/L}\Big (\prod_{\q_M^{} |\q} U_{\q_M^{}}\Big)
=\prod_{\q_L^{}|\q}U_{\q_L^{}}^{p^n}\ , $$
and then,
$$\{u\in E_L \mid \forall g\in
G,\lambda_n(g(u))=1\}=E_L\cap\dsp\prod_{\q_L^{}
|\q}U_{\q_L^{}}^{p^n}\supset \No_{M/L}(E_M)\ ,$$
so that
$$|H^0(\mathcal G,E_M)|=|E_L/\No_{M/L}(E_M)|\ge |E_L/\{u\in E_L\,
\vert \,\forall
g\in G,\lambda_n(g(u))=1\}|\ ,$$
and we finally have from the character theory side:
$$|H^0(\mathcal G,E_M)|\ge p^{n(|G|-1)-\delta}\ .$$


{\bf (5) Study of $ \hat H^1(\mathcal G,E_{M})$ and majoration of
$|I_{M}^\mathcal G/P_{M}^\mathcal G|$:}

Recal that $M:= FL$.
According to [4], chapter IX, \S1, since the cyclic extension $M/L$
is totally split at all infinite places, the Herbrand quotient
$q(\mathcal G,E_M)$ of the units is given by:
$$q(\mathcal
G,E_M) :=\frac{|\hat H^0(\mathcal
G,E_M)|}{|\hat H^1(\mathcal
G,E_M)|}=\frac{1}{[M:L]}=\frac{1}{p^{n+d}}\  ,$$
and this gives:
$$|\hat H^1(\mathcal G,E_M)|=p^{n+d}\,.\,|\hat H^0(\mathcal G,E_M)|\ge
p^{n+d}\,.\,p^{n(|G|-1)-\delta} = p^{n|G|+d-\delta}\ .$$
On the other hand, one has the canonical isomorphism:
$$
\hat H^1(\mathcal G,E_M) \simeq P_M^\mathcal G/P_L
$$
where $P_M^\mathcal G$ is the group of principal ideals of
$M$ which are invariant under $\mathcal G$.

So one obtains: $$
|P_M^\mathcal G/P_L|\ge p^{n|G|+d-\delta}.
$$
Thus, from the formula:\footnote{Since we have supposed that $\q$ is
totally split in $L/K$.}
$$
|I_M^\mathcal G/I_L|=\prod_{\q_L^{} \vert \q}e_{\q_L}=p^{n|G|}\ ,
$$
where $I_M$ is the group of fractional ideals of $M$, one deduces:
$$
|I_M^\mathcal G/P_M^\mathcal G|=\frac{|I_M^\mathcal
G/P_L|}{|P_M^\mathcal G/P_L|}=\frac{|I_M^\mathcal
G/I_L|\,.\,|I_L/P_L|}{|P_M^\mathcal
G/P_L|}=\frac{p^{n|G|}\,.\,|\cl_L|}{|P_M^\mathcal
G/P_L|}\le\frac{p^{n|G|}\,.\,|\cl_L|}{p^{n|G|+d-\delta}}\ ,
$$
that is:
$$|I_M^\mathcal G/P_M^\mathcal G|\le|\cl_L|\,.\,p^{\delta-d}\ .$$

Note that the number at the right hand side does not depend on $n$.


{\bf (6) Does the class $ c_L^{}$ principalize in $ M$ ?}

Here, we also suppose:

$\bullet \ \frob(\q_L^{}, H_L/L)=c_L^{}$.

$M'$ being the maximal subfield of $M/L$ in which $\q_L^{}$
totally splits, one has:
$$\q_L^{}=\prod_{\q'_M| \q_L^{} \text{ in }  M'/L}\q'_M;$$
and for all prime $\q'_M| \q_L^{}$ of $M'/L$, we have
$\q'_M=\q_M^{p^n}$,
where $\q_M^{}$ is the unique prime of $M$
dividing $\q'_M$.

Finally in $M$,
$$\q_L^{}=\Bigg(\prod_{\q'_M | \q_L^{}}\q_M^{}\Bigg)^{p^n}\ ,$$
with $\dsp\prod_{\q'_M | \q_L^{}}\q_M^{}\in I_M^\mathcal G$.
According to Section {\bf (5)}, one has:
$$|I_M^\mathcal G/P_M^\mathcal G|\le|\cl_L|\,.\,p^{\delta-d},$$
then:
$$
\Bigg(\prod_{\q'_M |\q_L^{}}
\q_M\Bigg)^{|\pcl_L|\,.\,p^{\delta-d}}\in P_M^\mathcal G\ .
$$\smallskip
Hence, in $\cl_M$, the extended class $c_M^{}$ of $c_L^{}$ satisfies both:
$$
c_M^{}=\overline\q_L^{}=\Bigg(\prod_{\q'_M |
\q_L^{}}\overline\q_M^{}\Bigg)^{p^n}\quad \text{ and } \quad
\Bigg(\prod_{\q'_M |
\q_L^{}}\overline\q_M^{}\Bigg)^{|\pcl_L|\,.\,p^{\delta-d}}=1\ .$$
So, $w$ being such that $|\cl_L|=p^w$, $c_L^{}$ is principal in $M$ under
the assumption:

\medskip
\centerline{$ n\ge w+\delta-d. $}

\noindent{\bf (7) Existence of $\q$:}

We just proved that $c_L^{}$ principalizes in $M$ when $n$ is large
enough
and when $\q$ (or $\q_L^{} | \q$) satisfies the following six
conditions:\vspace{-0.1cm}

$(1)$\  \ $\q\nmid p$,\vspace{-0.2cm}

$(2)$\  \ $\mu_{p^n}\subset K_\q^\times$,\vspace{-0.2cm}

$(3)$\  \ $E_K\subset U_\q^{p^n}$,\vspace{-0.2cm}

$(4)$\  \ $\q$ is totally split in $L/K$,\vspace{-0.2cm}

$(5)$\  \  $\frob(\q_L^{},
L[\mu_{p^n},E_L^{1/{p^n}}]/L)=\lambda_n$,\vspace{-0.2cm}

$(6)$\  \ $\frob(\q_L^{}, H_L/L)=c_L^{}.$\vspace{-0.2cm}

The definition of $\lambda_n\in\gal(L[\mu_{p^n},E_L^{1/p^n}]/L)$ shows
it is trivial on
$L[\mu_{p^n},E_K^{1/p^n}]$ then conditions (4) and (5) imply (2) and (3), so
we only study compatibility between conditions (1), (4), (5), (6). This
compatibility is possible if and only if $c_L^{}$ and $\lambda_n$ are
equal
on the extension $H_L\cap
L[\mu_{p^n},E_L^{1/p^n}]$ of $L$.

Let $m$ be the integer such that $|(\mu_L)_p|=p^m$; one has, where the
exponent $ab$ means abelian subextension over $L$:
$$H_L\cap
L[\mu_{p^n},E_L^{1/p^n}]\subset H_L\cap
(L[\mu_{p^n},E_L^{1/p^n}])^{ab}=H_L\cap
L[\mu_{p^{n+m}},E_L^{1/p^m}]\ ;$$
let $m'$ be the integer such that $H_L\cap
L[\mu_{p^\infty}]=L[\mu_{p^{m'}}]$, so $m'\ge m$ and
$$H_L\cap
L[\mu_{p^n},E_L^{1/p^n}]= H_L\cap
L[\mu_{p^{m'}},E_L^{1/p^m}]\ .$$
The exponent of the Galois group $\gal(L[\mu_{p^{m'}},E_L^{1/p^m}]/L)$
is less than $p^{m'}$, so is this of $\gal(H_L\cap
L[\mu_{p^n},E_L^{1/p^n}]/L)$. According to Section {\bf (1)}, taking $a=m'$,
one can suppose that $c_L^{}\in\cl_L^{p^{m'}}$ (replacing $L$ by $L'$ as
in Section{\bf (1)}; note that $m'(L')=m'(L)$ because $L'/L$ is unramified
at  all places
dividing $p$, so that $c_L^{}\in\cl_L^{p^{m'(L')}}$ as expected); in that
case, the restriction of $c_L^{}$ is trivial on $H_L\cap
L[\mu_{p^n},E_L^{1/p^n}]$.

About the restriction of $\lambda_n$ on $H_L\cap
L[\mu_{p^n},E_L^{1/p^n}]$, we can as well suppose it is trivial, by
replacing eventually $\lambda_n$ by $\lambda_n^{p^{m'}}$ ({\em i.e.}
$\varphi$ by $p^{m'}\varphi$) which has the same
properties.\footnote{See the Appendix.}

Up to replacing $L$ by $L'$ and choosing a convenient $\varphi$, the
restrictions of $c_L^{}$ and of $\lambda_n$ are both trivial on $H_L\cap
L[\mu_{p^n},E_L^{1/p^n}]$: \v{C}ebotarev theorem then ensures the
existence of infinitely many convenient primes $\q$ of $K$
satisfying  all
conditions, and each one gives us an abelian compositum $M$ in which
$L$ principalizes. This proves the Main Theorem.

\bigskip
\noindent{\large \bf III. Proofs of corollaries}

\smallskip
Corollary 1 is just the case $K=\mathbb Q$, and the Kronecker-Weber
theorem which states that abelian extensions of $\mathbb Q$ are
cyclotomic.

\smallskip
Corollary 2 is
equivalent to the following fact:
Let $K$ be a number field and $K^{ab}_+$ its maximal $\infty$-split
abelian extension. Any field $L$ containing
$K^{ab}_+$ and in which at least one infinite place totally splits
over  $K$ is principal.

To prove this, Let $\mathfrak a$ be a fractional ideal of finite type of
$L$. Of course $\mathfrak a$ is as well a fractional ideal of a subfield
$L_{\mathfrak a}$
of  $L$ with
finite dimension over $\mathbb Q$. We can suppose
$L_{\mathfrak a}\supset K$, then $L_{\mathfrak a}/K$ is an extension of
number fields in which
at least one infinite place is totally split. According to the main
theorem,  ${\mathfrak a}$
 principalizes in an abelian compositum $L_{\mathfrak a}F$ of $L_{\mathfrak
a}/K$; but 
$L_{\mathfrak a}F$ is contained in
$L_{\mathfrak a} K^{ab}_+ \subseteq L$ and so ${\mathfrak a}$ is principal
in $L$.

\smallskip
Corollary 3 comes from corollary 2, by taking $K=\mathbb Q$ and
$K=\mathbb  Q(i)$, respectively.

\smallskip
Corollary 4 comes from Corollary 2.

\bigskip
\centerline{D -- APPENDIX\,\footnote{Written  by G. Gras and J.-F. Jaulent.}}

\bigskip
The original project of publication of S. Bosca was first written in
french from his thesis and a provisional text, in english, was given to us
before his departure from the University.
Thus, due to the interest of the ideas of this work, it has been decided
to publish it, with suitable corrections in the text and with a complement
which is given below in this Appendix.

\bigskip
{\bf  Definition of $\varphi$.} For a group of global units
$E$ we denote by $\ov E$ the quotient of $E$ by its torsion subgroup
$\mu$.

Let $\No$ be the norm in $L/K$ and let ${}_{\No}\!E$ be the kernel of
$\No$ in $E$.
From the exact sequence:
$$1 \too {}_{\No}\!\ov E_L \tooo\ov  E_L \tooo \No (\ov E_L) \subseteq \ov
E_K \too 1 \ \
(\hbox {finite index}) $$
we get:
$$\mathbb Q \otimes_\mathbb Z (\ov  E_L) = \mathbb Q \otimes_\mathbb
Z ({}_{\No}\!\ov E_L)
\oplus \mathbb Q \otimes_\mathbb Z (\ov E_K)\  \hbox{\it  i.e. } \
\mathbb Q \otimes_\mathbb Z (\ov  E_L/\ov  E_K) = \mathbb Q \otimes_\mathbb
Z ({}_{\No}\!\ov E_L). $$
Since at least one infinite place of $K$ is totally split in the extension
$L/K$, the Dirichlet--Herbrand theorem implies that the character of
$\mathbb Q\otimes_{\mathbb Z} ({}_{\No}\!\ov E_L)$ contains
$\chi(\mathbb Q[G])-1$ and the representation
$\mathbb Q \oplus \mathbb Q \otimes_{\mathbb Z} ({}_{\No}\!\ov E_L)$
contains at least a representation $R$ isomorphic to
$\mathbb Q[G]$.

We can put $R = \mathbb Q  \otimes_{\mathbb Z} \langle\, \theta \,\rangle$
with $\theta = \rho\, .\, \varepsilon_*$ where $\rho \in \mathbb
Q^\times$,
$\rho \ne \pm1$, and where $\varepsilon_* \in {}_{\No}\!\ov E_L$ may be seen
as
a `` relative Minkowski unit ''.

Thus any element $u \in R$ is written, in a unique manner,
$u = \theta_{}^\omega$ with $\omega = \sum_{g\in G} \alpha_g\, g \in \mathbb
Q[G]$.
It follows that $u$ is a unit (in $\langle\, \varepsilon_*
\,\rangle_{\mathbb Z[G]}$)
if and only if $\sum_{g\in G} \alpha_g =
0$. We note that in this case the $\alpha_g$ can be taken in
$\frac{1}{m}\mathbb Z[G]$ for a suitable $m \in \mathbb Z$ (for
instance $m = \vert \, G\, \vert$, but if necessary we can adjust the
value of $m$ large enough; at the end of the reasoning, Bosca uses
this possibility); in any case $m$ depends only on
$L/K$ and not on $n$.

For the same reasons, the choice of $\varepsilon_*$ is not crucial
and $\langle\, \varepsilon_* \,\rangle_{\mathbb Z[G]}$ is not necessarily a
direct
summand in ${}_{\No}\!\ov E_L$.

The map $\varphi$ is then defined as follows:
noting that $$\ov E_L/{}_{\No}\!\ov E_L \,.\, \ov E_K =
\ov E_L/{}_{\No}\!\ov E_L \oplus \ov E_K$$ is killed by $\vert \, G\,
\vert$,
for $\varepsilon \in \ov E_L$, we have $\varepsilon_{}^{\vert  G \vert} =
\eta_* \,.\, \varepsilon_0$, $\eta_* \in {}_{\No}\!\ov E_L$,
$\varepsilon_0 \in \ov E_K$.

Working in $\mathbb Q \oplus \mathbb Q \otimes_\mathbb
Z ({}_{\No}\!\ov E_L)$, in which $R$ is a direct summand,
we associate with $\varepsilon \in \ov E_L$ the component of
$\varepsilon^m$ on $R$, of the form
$\theta_{}^\omega$, with $\omega = \sum_{g\in G} \alpha_g\, g \in \mathbb
Z[G]$, where $\sum_{g\in G} \alpha_g=0$, then we put:
$$\varphi(\varepsilon) = \alpha_1 \in \mathbb Z. $$
This map is trivial on $\ov E_K$ and defines an element of
$\hom(E_L/\mu_L\,.\,E_K,\mathbb Z)$ with the $G$-module action defined
as usual by:
$$\psi^h(x) := \psi(x^{h^{-1}}),\ \hbox{ for all } \psi \in
\hom(E_L/\mu_L\,.\,E_K,\mathbb Z)
\ \hbox{ and all }  h \in G.$$
It is  clear that $\varphi$ generates a
$\mathbb Z[G]$-submodule whose character is  $\chi(\mathbb Z[G])-1$.
More precisely, a straightforward computation gives
$\varphi^g(\varepsilon) = \alpha_g$ for any $g \in G$, thus
$\omega = \sum_{g\in G} \alpha_g\, g = \sum_{g\in G}
\varphi(\varepsilon^{g^{-1}})\, g$ .

In the sequel of the main text we will put $\omega := \phi(\varepsilon).$

At this step, Bosca introduces the map $\lambda_n$:
$$\lambda_n : \, E_L \tooo E_L/\mu_L\,.\,E_K \mathop{\tooo}^{\varphi}
 \mathbb Z \tooo \mu_{p^n}^{} \too 1$$
by a choice of a primitive $p^n$th root of unity.

This yields an element of $\hom(E_L ,\mu_{p^n}^{})$ which will be,
by abuse of notation,
 identified, via the Kummer duality between radicals and Galois groups,
to the corresponding element $\sigma'$ of the Galois group of
$L[\mu_{p^n},E_L^{1/p^n}]/L[\mu_{p^n}]$; then one creates a new
condition by saying that $\sigma'$ coincide with a suitable
Frobenius $\sigma$, which is the key idea for the proof of the conjecture
involving the necessary and sufficient condition about the splitting
of at least an infinite place.

\end{document}